\newif\ifpictures
\picturestrue
\newif\ifcomment
\commenttrue

\documentclass[12pt, reqno]{amsart}
\usepackage{latex_base}
\usepackage{macros} 

\usepackage{etoolbox}
\AtBeginEnvironment{theorem}{\setlist[enumerate, 1]{font=\upshape, nosep,  wide=0.5em, before=\leavevmode}}

\author{Sascha Timme} 

\address{Technische Universit\"at Berlin, Chair of Discrete Mathematics/Geometry, Germany\medskip} 
 
\email{timme@math.tu-berlin.de}


\title[Mixed Precision Path Tracking]{Mixed Precision Path Tracking for  \\Polynomial Homotopy Continuation}

\begin{document}

\begin{abstract}
This article develops a new predictor-corrector algorithm for numerical path tracking in the context of polynomial homotopy continuation.
In the corrector step it uses a newly developed Newton corrector algorithm which rejects an initial guess if it is not an approximate zero.
The algorithm also uses an adaptive step size control which builds on a local understanding of the region of convergence of Newton's method and the distance to the closest singularity.
To handle numerically challenging situations the algorithm uses mixed precision arithmetic.
The efficiency and robustness are demonstrated in several numerical examples.
\end{abstract} 

\maketitle

\section{Introduction}\label{sec:Introduction}

Systems of polynomial equations arise in many applications across the sciences including computer vision \cite{Hartley:Sturm:1997, Stewenius:Schaffalitzky:Niser:2005}, chemistry \cite{Morgan:2009}, kinematics \cite{Sommese:Wampler:2011} and biology \cite{Narula:Tiwari:Igoshin:2016}.
A numerical method for finding all isolated solutions of a system $F$ of $n$ polynomials in $n$ unknowns is homotopy continuation \cite{Sommese:Wampler:2005}.
This method constructs a homotopy $H(x,t): \C^n \times \C \rightarrow \C^n$ such that $H(x,t)$ is a polynomial system for all $t$, $H(x,1) = F(x)$ and $H(x,0)$ is a system whose isolated solutions are known.
There is a well-developed theory \cite{Sommese:Wampler:2005} on how to construct such homotopies to guarantee, with probability one,
that every isolated solution of $F(x) = 0$ is the endpoint of at least one smooth solution path
$x(t)$.
A solution path $x(t)$ is implicitly defined by the~conditions
\begin{equation}\label{eq:implicit_solution_path_introduction} 
H(x(t), t) = 0\; \text{ for all } t \in [0, 1) \text{ and } x(0) = x_0 \;.
\end{equation}
In order to trace a path $x(t)$, the problem \eqref{eq:implicit_solution_path_introduction} is treated as a sequence of problems
\begin{equation}\label{eq:seq_problems}
H(x(t_k), t_k) = 0\;,\quad k=0,1,2,\ldots 
\end{equation}
with an (a-priori unknown) subdivision $0 = t_0 < t_1 < \ldots < t_M = 1$ of the interval~$[0, 1]$.
Each of the problems \eqref{eq:seq_problems} is then solved by
a \emph{correction method}, usually Newton's method, under the assumption that a \emph{prediction method}, e.g., Euler's method, provides a good starting point.
The choice of step size $\Delta t_k = t_{k+1} - t_k$ is usually given by an \emph{adaptive step size control}.
The step size must be chosen appropriately:
if the step size is too large, the prediction can be outside the zone of convergence of the corrector, while a too small step size means progress is slow.
There have been many efforts to design such adaptive step size controls \cite{Schwetlick:Cleve:87,Kearfott:Xing:1994,Gervais:Sadiky:2004}.

In the context of polynomial homotopy continuation methods, two phenomena need particular attention.
Polynomial systems often have singular solutions, and thus, the paths leading to these solutions are necessarily ill-conditioned at the end.
While endgame methods \cite{Morgan:Sommese:Wampler:1990, Morgan:Sommese:Wampler:1992} exist to compute singular solutions, these still require to track the solution path sufficiently close to the singularity.
Usually, homotopies guarantee, with probability one, that no path passes through a singularity before reaching its endpoint. However, there is a non-negligible chance that a near-singular condition is encountered during the tracking.

Also, if two different solution paths are near to each other, then this can cause \textit{path jumping}.
That is, the solution that is tracked `jumps' from one path to another.
The typical reason is that starting from a point on the tracked path the prediction method returns a point which, according to the correction method, is a numerical approximation of a point that is on a different path than the tracked one.
A possible result of path jumping is that not all isolated solutions of a polynomial system are computed.

Therefore, path tracking algorithms are required to reduce the risk of path jumping
and they need to be able to handle ill-conditioned situations during the tracking.
Existing software packages, e.g., Bertini \cite{Bertini}, use a version of the following path tracking algorithm.
The algorithm has the following parameters: An initial step size $\Delta t >0$, a number of corrector iterations $N \ge 1$ allowed per step, a step adjustment factor $\lambda \in (0, 1)$, a step expansion integer $E \ge 1$, and a minimum step size $t_{\min}$.
Additionally, there is a tracking tolerance $\varepsilon > 0$. This means that for a given $t$ an approximate solution $x \approx x(t)$ has to satisfy a normwise absolute error $\norm{x - x(t)}_{\infty} \le \varepsilon$.

Given an approximate solution $x \approx x(t)$ the prediction method provides an initial guess $\hat{x}(\Delta t)\approx x(t+\Delta t)$. Then, Newton's method iteratively improves the approximation $\hat{x}(\Delta t)$.
If the required tracking tolerance $\varepsilon$ is achieved with a most $N$ iterations, then the solution is updated and $t = t + \Delta t$.
If there are $E$ successes in a row, then the step size is expanded to $\Delta t = \lambda^{-1}\Delta t$.
If on the other hand the tolerance is not achieved with at most $N$ iterations, then the step size is reduced to $\Delta t = \lambda \Delta t$. If $\Delta t < t_{\min}$, the algorithm terminates with a failure. Otherwise the procedure is repeated until $t=1$ is reached.

The key to avoid path jumping is to allow only a small number of Newton iterations, typically only $N=2$ or $N=3$. In practice, this is often sufficient for the initial guess $\hat{x}(\Delta t)$ to stay within a small enough region surrounding the path such that no path jumping occurs.
However, if two paths are closer than the required tracking tolerance $\varepsilon$ for some $t^* \in (0,1)$, then this algorithm tends to fail for these paths.
This is shown in the computational experiments in Section~\ref{sec:computations} for two different examples.
Therefore, it is necessary to choose the path tracking tolerance $\varepsilon$ smaller than the  minimal pairwise distance of any two paths.
However, knowing the optimal choice of $\varepsilon$ a-priori is impossible.
Thus, one has to use either a pessimistic value for $\varepsilon$ or resort to trial and error.
But choosing $\varepsilon$ small does not only slow down the tracking of \emph{all} paths, it also can result in new tracking failures.
The reason for this is that Newton's method in floating point arithmetic \emph{cannot} always produce solutions whose relative normwise error is smaller than $\varepsilon$.
This was shown by Tisseur in \cite{Tisseur:2001} and is explained in detail in Section~\ref{subsec:comp-newton}.

To avoid tracking failures due to insufficient precision Bates, Hauenstein, Sommese, and Wampler \cite{Bates:Hauenstein:Sommese:Wampler:2008} developed an adaptive precision version of the above described path tracking algorithm.
During the tracking, the algorithm dynamically changes the working precision such that Newton's method can theoretically always produce solutions accurate enough for the desired tracking tolerance.
This eliminates the problem of insufficient precision in exchange of a possibly high computational cost.
But it also still leaves open the problem of picking a suitable tolerance $\varepsilon$.

This article introduces a new robust path tracking algorithm, Algorithm \ref{alg:path-tracking}, which does not require the choice
of a path tracking tolerance $\varepsilon$ or a maximal number $N$ of corrector iterations allowed per step.
The key idea is to use a more intrinsic measure for accepting an initial guess in the Newton corrector:
An initial guess $\hat{x}(\Delta t)$ should only be accepted if the Newton iterates $\hat{x}(\Delta t) = \xn{0}, \xn{1}, \xn{2}, \ldots$ satisfy 
\begin{equation}\label{eq:intro-contraction}
    \norm{\xn{j+1} - \xn{j}} \le a ^{2^j - 1} \norm{\xn{1} - \xn{0}} 
\end{equation}
for $j=1,2,\ldots$ and some fixed constant $a \in (0, \frac12]$.
If the initial guess satisfies \eqref{eq:intro-contraction}, then it is an \emph{approximate zero}.
This notion was introduced by Smale \cite{Smale:1986} for $a = \frac12$ and plays an important role in the complexity analysis of polynomial homotopy continuation methods \cite{Buergisser:Cucker:2011,Lairez:2017}.
Based on this idea, this article develops a new Newton corrector algorithm, Algorithm \ref{alg:newton-corrector}. The algorithm rejects an initial guess if \eqref{eq:intro-contraction} is not satisfied for some $j=1,\ldots,m$ where $m$ is dynamically chosen as the maximal number of iterations for which \eqref{eq:intro-contraction} can be satisfied in fixed precision floating point arithmetic.

The proposed path tracking algorithm combines the new Newton corrector algorithm with an adaptive step size control that chooses $\Delta t$ based on local geometric informations.
The step size control extends an adaptive step size control developed by Deuflhard \cite{Deuflhard:79} and combines it with the insight of \cite{Telen:VanBarel:Verschelde:2019} to use Pad\'e approximants as prediction methods.
In particular, the algorithm builds a local understanding of the region of convergence of Newton's method and the distance to the closest singularity.
This makes the developed path tracking algorithm robust against path jumping.
To handle numerically challenging situations the algorithm uses \emph{mixed precision} arithmetic.
That is, while the bulk of the computations is performed in double precision some computations are performed, if necessary, in extended precision.
This article is accompanied by a prototype implementation of the algorithm and an implementation will also be available in version 2.0 of \texttt{HomotopyContinuation.jl} \cite{HomotopyContinuation.jl}.

This article is organized as follows. Section \ref{subsec:newton-convergence-results} reviews a Kantorovich style convergence theory of Newton's method and
Section \ref{subsec:comp-newton} develops a new Newton corrector algorithm, Algorithm \ref{alg:newton-corrector}, based on requirement \eqref{eq:intro-contraction}.
In Section \ref{sec:predictor} the use of Pad\'e approximants as prediction method is developed.
In Section \ref{sec:path-tracking-alg} the results from the previous sections are used to develop an adaptive step size control. Finally, the new path tracking algorithm, Algorithm \ref{alg:path-tracking}, is stated.
The algorithm's effectiveness and robustness is shown through several numerical
experiments in Section \ref{sec:computations}.

\section{Newton's Method: Theory and Computational Aspects}\label{sec:newton}
The path tracking algorithm for a path $x(t)$ consists of three main components:
An adaptive step size routine that provides a step size $\Delta t$, a predictor that produces an initial guess $\hat{x}$ of $x(t + \Delta t)$, and a corrector that takes $\hat{x}$ and returns either an approximation of $x(t + \Delta t)$ or rejects $\hat{x}$.
This section focuses on Newton's method as a corrector.
The goal is to understand the size of the region of convergence of Newton's method as well as the behavior of Newton's method in floating point arithmetic and to translate this into a Newton corrector algorithm.

\subsection{Convergence results}\label{subsec:newton-convergence-results}
Let $D\subseteq \C^n$ an open set. Let $F: D \subseteq \C^n \rightarrow \C^n$ be an analytic function and $J_F: \C^n \rightarrow \C^n$ its Jacobian. Consider the Newton iteration
\begin{equation}\label{eq:newton-iteration}
\begin{array}{rl}
J_F(\xn{j}) \Delta \xn{j} &= F(\xn{j}) \\
\xn{j+1} &= \xn{j} - \Delta \xn{j} \\ 
\end{array}
, \quad j=0,1,2,\ldots
\end{equation}
starting at the initial guess $\xn{0} \in D$.
In 1948 Kantorovich \cite{Kantorovich:1948} already showed sufficient conditions for the convergence of Newton's method and the existence
of solutions. He also showed the uniqueness region of solutions and provided error estimates.
A particular property of Newton's method is that the iterates \eqref{eq:newton-iteration} are invariant under general linear transformations of $F$.
That is, given a start value $\xn{0} \in D$ and $A \in {\mathrm{GL}}_n(\C)$ the Newton iterates of $AF(x)$ and $F(x)$ coincide. This property is referred to as \emph{affine covariance} \cite{Deuflhard:2011}.
In the following an affine covariant version of a Kantorovich style convergence theorem for Newton's method is given.
The statement is due to Deuflhard and Heindl \cite{Deuflhard:Heindl:79} with error bounds from Yamamoto \cite{Yamamoto:1985}.
\begin{theorem}[Newton-Kantorovich \cite{Deuflhard:Heindl:79,Yamamoto:1985}]\label{thm:newton}
Let $F: D \subseteq \C^n \rightarrow 
\C^n$ be analytic. For some $\xn{0} \in D$, assume that $J_F(\xn{0})$ is invertible and that for all $x,y \in D$
\begin{align*}
&\norm{J_F(\xn{0})^{-1} (J_F(x) - J_F(y)) } \le \omega \norm{x - y}, \\
&\norm{\Delta \xn{0}} = \norm{J_F(\xn{0})^{-1} F(\xn{0})} \le \beta
\end{align*}
and $h_0 := \omega \beta \le \frac12$.
Let $r^* = (1 - \sqrt{1 - 2h_0}) / \omega$ and $\bar{S}(\xn{0}, r^*) = \{ x \,|\, \norm{x - \xn{0}} \le r^* \} \subseteq D$.
Then:
\begin{enumerate}
\item The iterates \eqref{eq:newton-iteration} are well-defined, remain in $\bar{S}(\xn{0}, r^*)$ and converge to a solution \phantom{$x^*$ of} $x^*$ of $F(x)=0$.
\item The solution is unique in $S(\xn{0}, r^{**}) \cap D$ where $r^{**} = (1 + \sqrt{1 - 2h_0}) / \omega$.
\end{enumerate}
Furthermore, assume $h < \frac12$ and define the recursive sequence $h_j = \frac{h_{j-1}^2}{2(1 - h_{j-1})^2}$.
Then also the following error estimates hold.
\begin{equation}\label{eq:newton-iter-acc-bound}
\norm{\Delta \xn{j}} \le \frac12 \omega \frac{\sqrt{1 - 2h_j} }{\sqrt{{1 - 2h_0}}} \norm{\Delta \xn{j-1}}^2, \quad j=1,2,3,\ldots
\end{equation}
\begin{equation}\label{eq:newton-iter-acc-solution-bound}
\norm{\xn{j} - x^*} \le \frac{2\norm{\Delta \xn{j}}}{1 + \sqrt{1 - 2 \omega \frac{\sqrt{1 - 2h_j} }{\sqrt{{1 - 2h_0}}} \norm{\Delta \xn{j}}}}, \quad j=0,1,2,\ldots
\end{equation}
\end{theorem}

A drawback of the Newton-Kantorovich theorem is that it is not possible to obtain sufficient conditions for the convergence of Newton's method by only using data from the initial guess $\xn{0}$. Instead, local information about the Lipschitz constant $\omega$ is required.
The necessity of local information motivated Smale to develop his $\alpha$-theory \cite{Smale:1986}, which only requires data from the initial guess $\xn{0}$ to compute sufficient conditions for the convergence of Newton's method.
This point of view has valuable features for the theory of computation.
In particular, it is the building block for the complexity analysis of polynomial homotopy continuation methods \cite{Beltran:Pardo:2009,Buergisser:Cucker:2011,Beltran:Leykin:2013,Lairez:2017}, certified path tracking algorithms \cite{Beltran:Leykin:2013} and the a posteriori certification of zeros \cite{Hauenstein:Sottile:2012}.

In order to talk about initial guesses which are super-convergent starting with the first iteration, Smale \cite{Smale:1986} introduced the concept of an \emph{approximate zero}.
\begin{definition}[Approximate zero]\label{def:approx-zero}
The point $\xn{0} \in \C^n$ is an approximate zero of $F$ if the Newton iterate $\xn{j}$ is defined for $j = 1,2,\ldots$ and satisfies
$$
\norm{\Delta \xn{j}} \le \left( \frac12 \right)^{2^j - 1} \norm{\Delta \xn{0}}\,. 
$$
If $\xn{0}$ is an approximate zero, then the true zero $x^* \in \C^n$ of $F$ to which the iterates are converging is the \emph{associated zero} of $\xn{0}$.
\end{definition}
Smale's \(\alpha\)-theorem gives a sufficient condition for $\xn{0}$ to be an approximate zero.
The theorem uses
\begin{align*}
\gamma(F,x) &\,\,=\,\, \sup_{k\geq 2}\big\Vert \frac{1}{k!}\,J_F(x)^ {-1} D^kF(x)\big\Vert^\frac{1}{k-1}\ \text{ and } \beta(F,x) \,\,=\,\, \Vert J_F(x)^{-1}F(x)\Vert
\end{align*}
where $D^kF$ is the tensor of order-$k$ derivatives of \(F\) and the tensor $J_F^{-1}D^kF$ is understood as a multilinear map $A:(\mathbb{C}^n)^k\to \mathbb{C}^n$ with norm $\Vert A\Vert := \max_{\Vert v \Vert = 1} \Vert A(v,\ldots,v)\Vert$.
\begin{theorem}[Smale's $\alpha$-theorem \cite{Smale:1986}]\label{thm:alpha_theorem}
There is a naturally defined number $\alpha_0$ approximately equal to
0.130707 such that if $\alpha(F, \xn{0}) := \beta(F,\xn{0}) \, \gamma(F,\xn{0}) < \alpha_0$, then $\xn{0}$ is an approximate zero of $F$.
\end{theorem}
It is also possible to give sufficient conditions for $\xn{0}$ to be an approximate zero
under the assumptions of the Newton-Kantorovich Theorem \ref{thm:newton}.
\begin{lemma}\label{lemma:conv-speed}
Using notation from Theorem \ref{thm:newton} assume
$$h_0 = \omega \norm{\Delta \xn{0}} \le 2(\sqrt{4a^4 + a^2} - 2a^2) =: h(a) $$
for a parameter $0 < a < 1$. Then, the contraction factors
\begin{align}
\Theta_{j} := \frac{\norm{\Delta \xn{j+1}}}{\norm{\Delta \xn{j}}} &\le a^{2^{j}} \quad,\,j=0,1,2,\ldots \label{eq:contraction_factors}\\
\intertext{and the error bounds}
\norm{\Delta \xn{j}} &\le a^{2^j - 1} \norm{\Delta \xn{0}} \quad,\,j=1,2,\ldots \label{eq:conv2}
\end{align}
are satisfied.
In particular, $\xn{0}$ is an approximate zero if $h_0 \le \sqrt2 - 1$.
\end{lemma}
\begin{proof}
From the error estimate \eqref{eq:newton-iter-acc-bound} follows
\begin{align}
\Theta_j &\le \frac12 \omega \frac{\sqrt{1 - 2h_j} }{\sqrt{{1 - 2h_0}}} \norm{\Delta \xn{j}} \le \frac12 \omega \sqrt{\frac{1}{1-2h_0}} \norm{\Delta \xn{j}} 
\label{eq:theta_error_estimate} \\
&= \frac12 \omega \norm{\Delta \xn{0}} \sqrt{\frac{1}{1-2h_0}} \prod_{\ell=0}^{j-1}\Theta_j
= \frac{h_0}{2\sqrt{{1-2h_0}}} \prod_{\ell=0}^{j-1}\Theta_\ell \nonumber \,.
\end{align}
From $h_0 \le 2(\sqrt{4a^4 + a^2} - 2a^2) < \frac12$ follows $\frac{h_0}{2\sqrt{{1-2h_0}}} \le a$ and therefore
$$\Theta_j \le a \prod_{\ell=0}^{j-1}\Theta_\ell \le aa^{\sum_{\ell=0}^{j-1}2^\ell} = a^{2^j} \,.$$
Statement \eqref{eq:conv2} follows from \eqref{eq:contraction_factors} by observing
$$
\frac{\norm{\Delta \xn{j}}}{\norm{\Delta \xn{0}}} = \prod_{\ell=0}^{j-1} \Theta_\ell \le a^{2^j - 1}\,.
$$ 
\end{proof}

The Newton-Kantorovich theorem and Smale's $\alpha$-theorem both
give sufficient conditions for an initial guess to be an approximate zero.
For the Newton-Kantorovich theorem, a (local) estimate of the Lipschitz constant $\omega$ needs to be obtained and for Smale's $\alpha$-theorem $\gamma$ needs to be computed.
The path tracking algorithm developed in this paper is based on the Newton-Kantorovich theorem since a (rough) estimate of $\omega$ can be computed with almost no additional cost during the Newton iteration.

A computational estimate $[\omega]$ of $\omega$ is
\begin{equation}\label{eq:comp-omega}
[\omega] = 2\frac{\norm{\Delta \xn{1}}}{\;\norm{\Delta \xn{0}}^2} \,.
\end{equation}
This can be seen as follows.
Using the error estimate \eqref{eq:newton-iter-acc-bound}
$$
\norm{\Delta \xn{1}} \le \frac12 \omega \sqrt{\frac{{1 - 2h_1} }{{{1 - 2h_0}}}} \norm{\Delta \xn{0}}^2
$$
together with the observation
$\frac{1-2h_1}{1 - 2h_0} = {(1 - h_0)^{-2}}$
follows
\begin{equation*}
\norm{\Delta \xn{1}} \le \frac12 \omega\frac{1}{1 - \omega \norm{\Delta \xn{0}}} \norm{\Delta \xn{0}}^2
\end{equation*}
which is equivalent to
\begin{equation}\label{eq:comp-omega-help-ineq}
\frac{2\norm{\Delta \xn{1}}}
{\norm{\Delta \xn{0}}^2 + 2\norm{\Delta \xn{0}} \norm{\Delta \xn{1}} } \le \omega \,.
\end{equation}
The computational estimate \eqref{eq:comp-omega} is now obtained by upper bounding \eqref{eq:comp-omega-help-ineq} with
$$ \frac{2\norm{\Delta \xn{1}}}
{\norm{\Delta \xn{0}}^2 + 2\norm{\Delta \xn{0}} \norm{\Delta \xn{1}} } \le 2\frac{\norm{\Delta \xn{1}}} 
{\; \norm{\Delta \xn{0}}^2} = [\omega] \,.$$ 

\subsection{Computational Aspects and Floating Point Arithmetic}\label{subsec:comp-newton}
After establishing the theoretical foundations of Newton's method as well as a method to obtain a computational estimate of the Lipschitz constant $\omega$, these results are now used to guide the development of a Newton corrector algorithm.
For this, the behavior of Newton's method in floating-point arithmetic has to be taken into account.

\subsection*{Limit Accuracy} The following assumes the standard model of floating point arithmetic \cite[section 2.3]{Higham:2002}
$$
fl(x \; \mathrm{op} \; y) = (x \; \mathrm{op} \; y)(1 + \delta), \quad | \delta | \le u, \quad \mathrm{op}=+,-,*,/
$$
where $u$ is the unit roundoff. In standard double precision arithmetic $u = 2^{-53} \approx 2.2 \cdot 10^{-16}$.
In \cite{Tisseur:2001}, Tisseur analyzed the limit accuracy of Newton's method in floating-point arithmetic.
Let $x^* \in \C^n$ be a zero of $F$ with $J_F(x^*)$ non-singular, and let $\xn{0} \in \C^n$ be an approximate zero of $F$ with associated zero $x^*$.
In floating point arithmetic, we have
$$ \xn{j+1} = \xn{j} - (J_F(\xn{j}) + E_j)^{-1}(F(\xn{j}) + e_j) + \varepsilon_j $$
where
\begin{itemize}
\item $e_j$ is the error made when computing the residual $F(\xn{j})$,
\item $E_j$ is the error incurred in forming $J_F(\xn{j})$ and solving the linear system for $\Delta \xn{j}$,
\item $\varepsilon_j$ is the error made when adding the correction to $\xn{j}$.
\end{itemize}
Assume that $F(\xn{j})$ is computed in the possibly extended precision $\bar{u} \le u$ before rounding back to working precision $u$
and assume that there exists a function $\psi$ depending on $F$, $\xn{j}$, $u$ and $\bar{u}$ such that 
$$
\norm{e_j} \le u \norm{F(\xn{j})} + \psi(F, \xn{j}, u, \bar{u}) \,.
$$
Similarly, assume that the error $E_j$ satisfies
\begin{equation*}
\norm{E_j} \le u \phi(F, \xn{j}, n, u)
\end{equation*}
for some function $\phi$ that reflects both the instability of the linear solver and the error made when forming $J_F(\xn{j})$.
Then the following statement holds \cite[Corollary 2.3]{Tisseur:2001}.
\begin{theorem}[\cite{Tisseur:2001}]\label{thm:tisseur}
Let $\xn{0}$ be an approximate zero with associated zero $x^*$, $x^* \ne 0$, assume that $J_F(x^*)$ is non-singular, satisfies $u \, \kappa(J_F(x^*)) \le \frac18 \,$ and assume that for all $j$
$$
u \norm{J_F(\xn{j})^{-1}} \phi(F, \xn{j}, n, u) \le \frac{1}{8} \,.
$$
\noindent Then, Newton's method in floating point arithmetic generates a sequence of iterates $\xn{j+1}$ whose normwise relative error decreases until the first $j$ for which
\begin{equation}\label{eq:newton-accuracy-limit}
\frac{\norm{\xn{j+1} - x^* } }{ \norm{x^*}} \approx \frac{\norm{J_F(x^*)^{-1}}}{\norm{x^*}} \psi(F, x^*, u, \bar{u}) + u =: \mu(x^*, u, \bar{u}) \,.
\end{equation}
\end{theorem}
In the following the value $\mu(x^*, u, \bar{u})$ is referred to as the \emph{limit accuracy}.
Theorem \ref{thm:tisseur} shows that the limit accuracy is influenced by three factors: the working precision $u$, the accuracy of the evaluation of the residual (in possibly extended precision $\bar{u}$), and the conditioning of the Jacobian.
The essential consequence of this is that Newton's method \emph{cannot} always produce solutions whose normwise relative error is on the order of the working precision $u$.
From the error estimate \eqref{eq:newton-iter-acc-bound} follows that if for a given $j$
\begin{equation}\label{eq:termination-newton}
\norm{\Delta \xn{j}} \le \frac{\omega \norm{\Delta \xn{j-1}}^2}{2\sqrt{1-2h_0}} \le \mu(x^*, u, \bar{u}) \norm{x^*} 
\end{equation} then the normwise relative accuracy of $\xn{j}$ in floating point arithmetic is only of order $\mu(x^*, u, \bar{u})$.
Assume that for a given $j$ \eqref{eq:termination-newton} is satisfied. Then a computational estimate $[\mu]$ of $\mu(x^*, u, \bar{u})$
can be obtained by computing $\norm{\Delta \xn{j}} / \norm{\xn{j+1}}$.

The following lemma shows that using extended precision improves the limit accuracy.
\begin{lemma}
    For extended precision $\bar{u} \le u$ it holds $\mu(x^*, u, \bar{u}) \approx \mu(x^*, u, u) \frac{\bar{u}}{u} + u$.
\end{lemma}
\begin{proof}
    From Section 4.3.2 in \cite{Bates:Hauenstein:Sommese:Wampler:2008} follows that for a system of polynomials given as a straight line program $\psi(F, \xn{j}, u, \bar{u})$ in Theorem \ref{thm:tisseur} is a linear function in $\bar{u}$, i.e., $\psi(F, \xn{j}, u, \bar{u}) = \frac{\bar{u}}{u} \psi(F, \xn{j}, u, u)$.
    The statement then follows from \eqref{eq:newton-accuracy-limit}.
\end{proof}
If the working precision $u$ is standard double-precision arithmetic, then computing with extended precision can be accomplished by using double-double arithmetic.
A double-double number is represented as an unevaluated sum of a leading double and a trailing double, resulting in a unit roundoff of $2^{-106} = u^2$.
Bailey \cite{Bailey:1995} pioneered double-double arithmetic, and implementations are nowadays available for a wide variety of programming languages and architectures.

Assume that for a fixed parameter $a \in (0,\frac12]$ the Newton iterates starting at the initial guess $\xn{0}$ are required to satisfy the contraction factors
\begin{equation*} 
\Theta_{j} = \frac{\norm{\Delta \xn{j+1}}}{\norm{\Delta \xn{j}}} \le a^{2^{j}} \; j=0,1,2,\ldots\,.
\end{equation*}
If the Newton iterates are computed with precision $\bar{u} = u$
then \eqref{eq:termination-newton} implies together with Lemma \ref{lemma:conv-speed} that if 
\begin{equation}\label{eq:newton-min-acc}
\omega \mu(x^*, u, u) > a^{2^{k} - 1} h(a) \norm{x^*}
\end{equation}
then there does not need to exist an initial guess $\xn{0}$ such that the first $k$ contraction factors are satisfied.
Given a fixed $k$, for instance, $k=2$, this gives a criterion when to use extended precision.
Similarly, if the Newton iteration is performed with extended precision, then it is possible to use only working precision again if 
$\omega \mu(x^*, u, u) < a^{2^{k} - 1} h(a)$.

The working precision $u$ is insufficient if the combination of the error in the evaluation of the Jacobian and the instability in the linear system solver become too large.
In this case, a multi-precision path tracking algorithm as \cite{Bates:Hauenstein:Sommese:Wampler:2008} is necessary.
However, as demonstrated in Section \ref{sec:computations}, using only double precision arithmetic for the linear system solver is sufficient for most applications.
Nevertheless, even if the precision $u$ is sufficient to achieve the limit accuracy, the analysis of Tisseur also shows that the convergence speed of Newton's method can decrease due to a too unstable linear system solver.
In this case, the theoretical convergence speed may not be achieved which in turn can lead to not satisfying the required contraction factors.
To circumvent this, the Newton updates are improved using mixed precision iterative refinement \cite{Higham:1997} if $\bar{u} < u$.
This stabilizes the linear system solver sufficiently to achieve the theoretical convergence speed.

\subsection*{Stopping Criteria}
Criteria for stopping the Newton iteration are now derived.
Assume that $\omega$ and the limit accuracy $\mu = \mu(x^*,u,\bar{u})$ are known.
If for any $j$ the contraction factor
\begin{equation}\label{eq:newt-rejection}
\Theta_{j} = \frac{\norm{\Delta \xn{j+1}}}{\norm{\Delta \xn{j}}} \le a^{2^{j}}
\end{equation}
is not satisfied then the iteration is stopped and the initial guess is rejected.
The iteration is stopped successfully at step $j$
if the next update would be smaller than the limit accuracy, i.e., 
\begin{equation}\label{eq:newt-acceptance}
    \frac{\omega \norm{\Delta \xn{j-1}}^2}{2\sqrt{1-2h(a)}} \le \norm{\xn{j}} \mu \,.
\end{equation}
An additional Newton update is computed to obtain a computational estimate of $\mu$.

\subsection*{Scaling} Before the full Newton corrector algorithm is stated, a final point is addressed.
So far, a simple rescaling of variables can change the behavior of the algorithm since $\omega$, $\mu$ and $\norm{\Delta 
\xn{j}}$ are not invariant under rescaling of variables.
Additionally, if $x^* = 0$ the accuracy needs to be measured with an absolute, and not a relative, normwise error.
A rescaling of variables is formally the change of variables
$$y = D^{-1}x, \quad D = \mathrm{diag}(d_1, \ldots, d_n), \, d_i \in \R_{>0} \,.$$
With $x = (x_1, \ldots, x_n)$ and $|x_i| \ne 0$ the choice $d_i \approx |x_i|$ results in new coordinates $y_i$ of unit order.
To deal with the case $|x_i| = 0$ as well as with possible overflows in floating point arithmetic an absolute threshold value $d_{\min} > 0$ of the form
\begin{equation}\label{eq:weight-init}
d_i = \max \{ |x_i|, d_{\min} \}
\end{equation}
has to be imposed. For instance, $d_{\min} = \max ( \sqrt{u} \max_i|x_i|, u)$.
To not introduce rounding errors, the scaling factors $D$ should be powers of the floating-point radix $\beta$ ($\beta= 2$ in the case of IEEE-754 floating point standard arithmetic). 
Instead of performing the change of variables explicitly in Newton's method the size of the Newton updates can also be measured with the \emph{weighted error}
$ \norm{D^{-1} \Delta \xn{j} } \,.$
Using the scaling factors $D$ allows the algorithm to perform independent of the initial provided variable scaling (assuming that the initial scaling is not too extreme).

\subsection*{The Algorithm}
Finally, a new Newton corrector algorithm, Algorithm \ref{alg:newton-corrector}, is stated.
The algorithm builds on the results developed in this section.
The idea of the algorithm is to reject an initial guess $\xn{0}$ if the Newton iterates $\xn{0}, \xn{1}, \xn{2}, \ldots$ do not satisfy 
\begin{equation*}
    \norm{\xn{j+1} - \xn{j}} \le a ^{2^j - 1} \norm{\xn{1} - \xn{0}} 
\end{equation*}
for $j=1,2,\ldots,m$ and some fixed constant $a \in (0, \frac12]$.
Here, $m$ is decided dynamically based on equation \eqref{eq:newt-acceptance}.
The rejection of an initial guess is performed using the slightly stricter criterion \eqref{eq:newt-rejection}.
The algorithm needs as input estimates of the limit accuracy $\mu$ and the Lipschitz constant $\omega$ and also returns updated estimates of these quantities.
During the path tracking, estimates are available by using the returned estimates from the previous steps. What to do at the beginning of the tracking, if these estimates are not available, is addressed after the algorithm.

\stepcounter{environment}
\begin{algorithm}[ht]
    \small
\caption{Newton Corrector}\label{alg:newton-corrector}
\begin{algorithmic}[1]
\Require $F: \C^n \rightarrow \C^n$, $\xn{0} \in \C^n$, $a \in (0,\frac12]$, estimate $[\mu] > 0$ of the limit accuracy $\mu$, estimate $[\omega] > 0$ of $\omega$, evaluation precision $\bar{u}\le u$, and positive scaling factors $D$ such that the coordinates of $D^{-1}\xn{0}$ are of unit order.
\Ensure Boolean indicating whether $\xn{0}$ was accepted, approximation $\bar{x} \in \C^n$ of a zero $x^*$ of $F$, updated estimate of the (limit) accuracy $\mu$ at $x^*$, updated estimate of $\omega$, number of updates $j$ and last contraction factor $\Theta_{j-2}$.
\Procedure{Newton}{$F, \xn{0}, a, [\mu], [\omega], \bar{u}, D$}
\State $j \gets 0$
\While{\textbf{true}}
\State $r \gets $ Evaluate $F$ at $\xn{j}$ with precision $\bar{u}$ and round result to precision $u$
\State Solve $J_F(\xn{j}) \Delta \xn{j} = r$ 
\State $\xn{j+1} \gets \xn{j} - \Delta \xn{j}$
\If{$j = 1$}
\State $[\omega] \gets 2\frac{\norm{D^{-1} \Delta \xn{1}}}
{\norm{D^{-1}\Delta \xn{0}}^2 } $
\Comment{Compute $\omega$ estimate} 
\EndIf
\If{$j \ge 1$ \textbf{and} $ \frac{\norm{D^{-1}\Delta \xn{j}}}{\norm{D^{-1} \Delta \xn{j - 1}}} > a^{2^{j - 1}}$}
\Comment{Check sufficient contraction}
\State \textbf{return} $(\textbf{false}, \xn{j+1}, [\mu], [\omega], j + 1,
\frac{\norm{D^{-1} \Delta \xn{j}}}{\norm{D^{-1} \Delta \xn{j - 1}}})$
\ElsIf{$\frac{\omega \norm{D^{-1} \Delta \xn{j}}^2}{2\sqrt{1-2h(a)}} \le [\mu]$} \Comment{Approaching limit accuracy}
\State $r \gets $ Evaluate $F$ at $\xn{j+1}$ with precision $\bar{u}$ and round result to precision $u$
\State Solve $ J_F(\xn{j+1}) \Delta \xn{j+1} = r$ 
\State $\xn{j+2} \gets \xn{j+1} - \Delta \xn{j+1}$
\State $[\mu] \gets \norm{D^{-1} \Delta \xn{j+2}}$ \Comment{Update of the limit accuracy}
\State \textbf{return} (\textbf{true}, $\xn{j+2}$, $[\mu]$, $[\omega]$, j+2, $\frac{\norm{D^{-1} \Delta \xn{j+1}}}{\norm{D^{-1} \Delta \xn{j}}}$)
\EndIf
\State $j \gets j + 1$
\EndWhile
\EndProcedure
\end{algorithmic}
\end{algorithm}
The algorithm requires estimates of the limit accuracy $\mu$ and the Lipschitz constant $\omega$.
During the path tracking, the computational estimates for both of them are available by using the computed estimates of the Newton corrector from the previous step.
However, this leaves open what to do for the first step.
There are two possibilities.
If the start solution is the solution of a previous tracking, then computational estimates of $\mu$ and $\omega$ are already available.
If this is not the case, the following heuristic, Algorithm \ref{alg:model-initialization}, to determine values for $[\mu]$ and $[\omega]$ proved to be helpful.
The idea is to add a small perturbation to the provided start solution and to perform two Newton steps. If the perturbation is sufficiently small, then the perturbed solution still converges to the provided start solution, and an estimate of $[\omega]$ and $[\mu]$ can be obtained.
As an added benefit, this provides a test to point out invalid start solutions, e.g., due to user error. For simplicity it is assumed that it is sufficient to compute the residual with precision $u$.
\stepcounter{environment}
\begin{algorithm}
\footnotesize
\caption{Model Initialization Heuristic}\label{alg:model-initialization}
\begin{algorithmic}[1]
\Require Candidate $\xn{0} \in \C^n$, $a \in (0,1)$, scaling factors $D$ such that the coordinates of $D^{-1}\xn{0}$ are of unit order.
\Ensure Boolean indicating whether the initialization was successful, estimate $[\mu]$ of the limit accuracy $\mu$ of the associated zero of $\xn{0}$, and an estimate $[\omega]$ of $\omega$.
\Procedure{ModelInitialization}{$\xn{0}$, $a$, $D$}
\State $v \gets \norm{D^{-1} J_F(\xn{0})^{-1} F(\xn{0})} + u$
\State $\varepsilon \gets \sqrt{v}$
\For{$k \gets 1:3$} \Comment{Try up to 3 different sizes of perturbations}
\State $\bar{x} \gets \xn{0} + \varepsilon D$
\Comment{Add relative perturbation} 
\State $\Delta_0 \gets J_F(\bar{x})^{-1} F(\bar{x})$
\State $\xn{1} \gets \bar{x} - \Delta_0$
\State $\Delta_1 \gets J_F(\xn{1})^{-1} F(\xn{1})$
\State $\xn{2} \gets \xn{1} - \Delta_1$
\If{$\norm{D^{-1} \Delta_1} / \norm{D^{-1} \Delta_0} < a$}
\State $[\omega] \gets 2\frac{D^{-1} \norm{\Delta \xn{1}}}
{\;\,\norm{D^{-1} \Delta \xn{0}}^2 } $
\State $[\mu] \gets \norm{D^{-1} \Delta_1}$
\State \textbf{return} (\textbf{true}, $[\mu]$, $[\omega]$)
\Else
\Comment{$\omega$ larger than $\varepsilon$, reduce $\varepsilon$}
\State $\varepsilon \gets \varepsilon u^{2^{-k}}$
\EndIf 
\EndFor
\State \textbf{return} (\textbf{false}, $[\mu]$, $[\omega]$)
\EndProcedure
\end{algorithmic}
\end{algorithm}

\section{Predictors and Pad\'e Approximants}\label{sec:predictor}
After carefully studying the Newton corrector in the previous section, the attention now shifts to the predictor.
Recall that the role of the predictor is to produce for a given step size $\Delta t$ an initial guess $\xn{0}$ such that the corrector converges sufficiently fast.
The choice of the step size $\Delta t$ will be addressed in the next section, but before it is essential to understand the influence of $\Delta t$ on the distance of the initial guess $\xn{0}$ to the solution path.

Consider the homotopy $H(x,t): \C^n \times \C \rightarrow \C^n$ and a constant $\bar{t} > 0$.
Given a solution $s \in \C^n$ of the system $H(x,0)=0$ assume that there is a solution path $x(t): [0,\bar{t}\,] \rightarrow \C^n$ implicitly defined by the conditions
\begin{equation}\label{eq:implicit_solution_path}
H(x(t), t) = 0\; \text{ for all } t \in [0, \bar{t}\,] \text{ and } x(0) = s \;.
\end{equation}
Also assume that $H_x(x(t),t)$ is nonsingular for all $t \in [0,\bar{t}\,]$.
Then $x(t)$ can be extended to a holomorphic function with $H(x(t^*), t^*) = 0$ for all $t^*$ in some nonempty open neighborhood of 0.

Without loss of generality, in the following, only the situation at $t=0$ is considered, thus $\Delta t = t$.
A predictor generates a \emph{prediction path} $\hat{x}(t): [0,\bar{t}\,] \rightarrow \C^n$ with $\hat{x}(0) = x(0)$ and they can be classified by the local order of the prediction error $\norm{\hat{x}(t) - x(t)}$.
\begin{definition}[Local order of a predictor]\label{def:order}
A predictor is of local order $p$ if there exists a $\tau > 0$ and a constant $\eta_p \ge 0$ such that for all $t \in [0,\tau]$
\begin{equation*}\label{eq:def_order_predictor}
\norm{\hat{x}(t) - x(t)} \le \eta_p t^p\,.
\end{equation*}
The constant $\tau$ is the \emph{trust region} of the predictor.
\end{definition}
\begin{example}
The Euler predictor $\hat{x}(t) = x(0) + t \dot{x}(0)$ is of order $p=2$ with $\tau = \bar{t}$ since
\begin{equation*}
\norm{x(t) - \hat{x}(t)} =
\norm{x(t) - x(0) - t \dot{x}(0)} 
\le \frac12 \max_{t \in [0, \bar{t}\,]} \norm{ \ddot{x}(t)} t^2 \,.
\end{equation*}
\end{example}

There are many different families of predictors known in the literature, with the most famous ones probably being (embedded) Runge-Kutta methods.
In \cite{Bates:Hauenstein:Sommese:2011} it is shown that for polynomial homotopy continuation higher-order Runge-Kutta methods are substantially more efficient than the Euler predictor.
However, in the following another particular class of predictors is considered: Pad\'e approximants.
See \cite{Baker:Graves-Morris:1996} for an exhaustive treatment of Pad\'e approximants.

\begin{definition}[Pad\'e approximant] \label{def:pade}
Let $x(t) = \sum_{\ell=0}^\infty c_\ell t^\ell$ be a convergent power series. The \textup{type $(L,M)$ Pad\'e approximant} is the rational function
\begin{equation*}
[L/M]_x = \frac{
a_0 + a_1 t + a_2 t^2 + \cdots + a_L t^L
}{
1 + b_1 t + b_2 t^2 + \cdots + b_M t^M
}
\end{equation*}
such that $x(t)$ and $[L/M]_x $ (considered as formal power series) satisfy
\begin{equation*}
[L/M]_x - x(t) \in \mathcal{O}(t^{L+M+1}) \;.
\end{equation*}
\end{definition}
\begin{remark}
A type $(L,M)$ Pad\'e approximant is a predictor of order $L+M+1$.
\end{remark}

\cite{Schwetlick:Cleve:87} demonstrates the effectiveness of a $(2,1)$ Pad\'e approximant as a predictor in path tracking algorithms. 
However, Pad\'e approximants are not only efficient predictors. They have the additional property to indicate the distance to the most nearby singularity.
This property was recently highlighted by Telen, Van Barel, and Verschelde in \cite{Telen:VanBarel:Verschelde:2019} where it is used to develop a path tracking algorithm that is robust against path jumping.

Since $x(t)$ is holomorphic in a nonempty open neighborhood of 0, there is a coordinatewise expansion of $x(t)$ as a convergent power series around 0.
Write $x_j(t) = \sum_{\ell = 0}^\infty c_\ell t^\ell$ for the Taylor expansion of the coordinate function $x_j(t)$ at $0$.
For sufficiently large $L$ the pole of the Pad\'e approximant $[L/1]_{x_j}$ indicates the \emph{distance} to the nearest singularity (also if it is a branch point). This is seen as follows.
A computation shows that if $c_L \neq 0$,
\begin{equation}\label{eq:pade-L}
[L/1]_{x_j} = c_0 + c_1 t + \ldots + c_{L-1} t^{L-1} + \frac{c_L t^L}{1-tc_{L+1}/c_L}.
\end{equation}
Hence the pole of $[L/1]_{x_j}$ is $c_L/c_{L+1}$ (or it is $\infty$ if $c_{L+1} = 0$).
Fabry's ratio theorem \cite{Fabry:1896} now states that if the limit $\lim_{L \rightarrow \infty} c_L / c_{L+1}$ exists it is a singularity of $x(t)$.
\begin{theorem}
Suppose that the coefficients of the power series $\sum_{\ell=0}^{\infty}c_\ell t^\ell$ are such that the limit $\lim_{L \rightarrow \infty} c_L / c_{L+1} = \lambda \ne 0$ exists.
Then, the series converges uniformly inside the disk $\{ |t| < |\lambda| \}$
and $\lambda$ is a singular point of $x_j(t) = \sum_{\ell=0}^{\infty}c_\ell t^\ell$.
\end{theorem}
For a fixed $L$ the modulus $|c_L / c_{L+1}|$ therefore can be assumed to be an approximation of the distance to the nearest singularity of $x(t)$ and a computational estimate of the trust-region $\tau$ of the Pad\'e approximant.
For a more extensive treatment of Pad\'e approximants in the context of homotopy continuation, see Section 3 of \cite{Telen:VanBarel:Verschelde:2019}.

For the computation of a Pad\'e approximant of type $(L,M)$ it is necessary to compute the local derivatives $x^{(\ell)}(0)$ for $\ell=1,\ldots,L+M$. For this Mackens \cite{Mackens:89} observed the following useful identity.
\begin{lemma}[\cite{Mackens:89}]
    The local derivatives $x^{(\ell)}(t)$ can be computed using the formula
    $$x^{(\ell)}(t) = - H_x(x(t), t)^{-1} R_\ell(t) \,,$$
    where 
    \begin{equation}\label{eq:higher-derivatives-identity}
    \begin{array}{rl}
    R_\ell(t) &= \left( 
    \frac{d}{d\lambda}
    \right)^\ell
    H \restr {\left(
    x(t) + \sum_{i=1}^{\ell - 1} \frac{1}{i!}x^{(i)}(t) \lambda^i,
    t + \lambda 
    \right)}{{\lambda=0}}\,.
    \end{array}
    \end{equation}
\end{lemma}
In \cite{Mackens:89} this identity is used for the computation of $x^{(\ell)}(0)$ by numerical differentiation.
A downside of numerical differentiation is that it can suffer from catastrophic cancellation resulting in useless results.
Instead of using numerical differentiation the expression \eqref{eq:higher-derivatives-identity} can be computed efficiently and accurately by using automatic differentiation \cite[Chapter 13]{Griewank:2008}.
In particular, the cost of computing $R_\ell$ using automatic differentiation is at most
$2\ell^2 + \mathcal{O}(\ell)$ times the cost of evaluating $H$ by a straight-line-program.
The dominating factor for the accuracy of $x^{(\ell)}(t)$ is the forward error of the linear system solving.
To ensure that computed derivatives are sufficiently accurate the forward error of the linear system solving should be monitored and if necessary be reduced by
using mixed precision iterative refinement \cite{Higham:1997}.

A robust Pad\'e approximant implementation also needs to handle the edge cases that $x_j^{(\ell)}(0) = 0$ for some $1\le j \le n$ and $1 \le \ell \le L+M$.
In \cite{Gonnet:Guttel:Trefethen:2013} a robust algorithm is proposed for computing Pad\'e approximants.
The provided implementation uses this algorithm for the computation of the Pad\'e approximants.

It is also possible to obtain an estimate of the local approximation error of a Pad\'e approximant. 
By comparing for each coordinate function $x_j(t)$ the coefficient of $t^{L+M+1}$ in
$$(a_0 + a_1 t + \ldots + a_L t^L) - (1 + b_1 t + \ldots + b_M t^M)(c_0 + c_1t + c_2t^2 + \ldots)$$
it follows
\begin{equation*}
e_{0,j} = - (c_{L+M+1} + b_1 c_{L+M} + \ldots + b_M c_{L+1}) \,.
\end{equation*}
Considering the Taylor expansion of $[L/M]_{x_j}$ at $0$ it follows
$$
x_j(t) - [L/M]_{x_j}(t) = e_{0,j}t^{L+M+1} + \mathcal{O}(t^{L+M+2}) \,.
$$
Therefore for a Pad\'e approximant a computational estimate $[\eta_{L+M+1}]$ of $\eta_{L+M+1}$ is
\begin{equation}\label{eq:errcoeff}
[\eta_{L+M+1}] = \norm{D^{-1}(e_{0,1}, \ldots, e_{0,n})}
\end{equation}
where $D > 0$ are the same scaling factors as used for the Newton corrector in Section \ref{sec:newton}.

\section{Step Size Control and Path tracking algorithm}\label{sec:path-tracking-alg}
After studying Newton's method as a correction method in Section \ref{sec:newton} and Pad\'e approximants as a prediction method in Section \ref{sec:predictor}, the results are now combined to derive an adaptive step size control. Afterward, the path tracking algorithm is stated.

Consider the homotopy $H(x,t): \C^n \times \C \rightarrow \C^n$ with a given solution $s \in \C^n$ of the system $H(x,0)=0$ and constant $\bar{t} > 0$.
Assume there is a solution path $x(t): [0,\bar{t}] \rightarrow \C^n$ implicitly defined by the conditions $x(0) = s$ and \eqref{eq:implicit_solution_path}, and assume that $H_x(x(t),t)$ is non-singular for all $t \in [0, \bar{t}\,]$. Denote by $\hat{x}(t): [0,\bar{\tau} \,] \rightarrow \C^n$ the prediction path produced by the Pad\'e approximant.
As in Section \ref{sec:predictor} only the situation at $t=0$ is considered such that $\Delta t = t$.

The goal of the step size routine is to provide a step size $t$ such that the Newton iterates $\xn{j}$ starting at the initial guess $\hat{x}(t) = \xn{0}$ satisfy for $j=0,1,2,\ldots$ the contraction factors
\begin{equation}\label{eq:contr-factor-step-size}
\Theta_{j} = \frac{\norm{\Delta \xn{j+1}}}{\norm{\Delta \xn{j}}} \le a^{2^{j}}
\end{equation}
for a fixed parameter $a \in (0, 1)$, for instance $a = 0.2$.
Recall from Lemma \ref{lemma:conv-speed} that if $a \le \frac12$ then $\hat{x}(t)$ is an approximate zero.
Assuming knowledge about the Lipschitz constant $\omega$ in a neighborhood of the path $x(t)$ and the theoretical quantities introduced in Section \ref{sec:predictor} it is possible to give a maximal theoretical feasible step size $t_{\max}$ such that this is the case.
The approach to use the theoretical quantities $\omega$, $\eta_p$, and $\tau$ to determine a maximal feasible step size such that Newton's method converges was pioneered by Deuflhard in \cite{Deuflhard:79}.
\begin{theorem}[\cite{Deuflhard:Heindl:79}]\label{thm:step-size}
Let $D \subseteq \C^n$ such that for all $x \in D$ and $t \in [0,\bar{t}\,]$, $\bar{t} > 0$, the Jacobian $H_x(x, t)$ is non-singular.
Assume that for each $t \in [0,\bar{t}\,]$ there exists a convex subset $D(t) \subseteq D$ with $x(t) \in D(t)$ where $x(t)$ denotes the unique solution path in $D \times [0,\bar{t}\,]$.
Let $\hat{x}(t): [0, \bar{t}\,] \rightarrow D$ denote a prediction path of order $p$ with trust-region $\tau$, i.e., with
$$ \norm{\hat{x}(t) - x(t)} \le \eta_p t^p \; \text{ for all $t \in [0,\tau]$} \,.$$
Moreover, assume for all $t \in [0,\bar{t}\,]$ the affine covariant Lipschitz condition
$$
\norm{ H_x(\hat{x}(t),t)^{-1} ( H_x(u,t) - H_x(v,t) ) } \le \omega \norm{ u - v}\text{ for all } u,v \in D(t)\,.
$$
For fixed $h \le \frac12$ let $t_{\max} = \min (t^*, \tau, \bar{t})$ where
\begin{equation}\label{eq:feasible-step-size}
t^* := 
\left(
\frac{\sqrt{1 + 2h} - 1}{ \omega \eta_p } 
\right)^{1/p}
\end{equation}
and for all $t \le t_{\max}$ let $B(t)$ denote a ball around $\hat{x}(t)$ with radius $(1 - \sqrt{1 - 2h}) / \omega$ and assume $B(t) \subseteq D(t)$.

Then, for all step sizes $t \le t_{\max}$ the Newton iterates starting at $\hat{x}(t) = \xn{0}$ are well-defined, remain in $B(t)$, converge towards ${x}(t)$ and satisfy $\norm{H_x(\hat{x}(t), t)^{-1} H(\hat{x}(t),t)} \le \frac{h}{\omega}$.
\end{theorem}
\begin{proof}
    For $h = \frac12$ the statement is Theorem 1.3 in \cite{Deuflhard:79}.
    The more general maximal step size \eqref{eq:feasible-step-size} and the inequality $\norm{H_x(\hat{x}(t), t)^{-1} H(\hat{x}(t),t)} \le \frac{h}{\omega}$ follows from equations (1.14a) and (1.14b) in the proof of Theorem 1.3 in \cite{Deuflhard:79}.
\end{proof}
In Theorem \ref{thm:step-size} there is a choice of the parameter $h \le \frac 12$. If this is sufficiently small then the contraction factors \eqref{eq:contr-factor-step-size} are satisfied. The following corollary makes this precise.

\begin{corollary}\label{cor:step-size}
In Theorem \ref{thm:step-size} choose $h \le h(a) = 2(\sqrt{4a^4 + a^2} - 2a^2)$ for $a \in (0,1)$. Then, the Newton iterates starting at $\hat{x}(t)$ satisfy the contraction factors
$$\frac{\norm{\Delta \xn{j+1}}}{\norm{\Delta \xn{j}}} \le a^{2^{j}} \,.$$
\end{corollary}
\begin{proof}
Since $\norm{H_x(\hat{x}(t), t)^{-1} H(\hat{x}(t),t)} \le \frac{h}{\omega}$ follows $h_0 \le h$ and using Lemma \ref{lemma:conv-speed}
follows the statement.
\end{proof}

If the step size is chosen according to Theorem \ref{thm:step-size}, then path jumping cannot happen.
However, to obtain the theoretical quantities $\omega$, $\eta_p$ and $\tau$ is very hard.
Instead the theoretical quantities are replaced by easy to obtain computational estimates $[\omega]$, $[\eta_p]$ and $[\tau]$.
Using the computational estimates and Corollary \ref{cor:step-size}, an estimate $[t_{\max}]$ of the maximal feasible step size $t_{\max}$
is given by
\begin{equation}\label{eq:comp-max-step-size}
[t_{\max}] = \min \left( \, [t^*], \, \beta_\tau [\tau], \, \bar{t} \, \right) 
\end{equation}
where
\begin{equation*}\label{eq:comp-max-feasible-partial}
[t^*] := 
\left(
\frac{\sqrt{1 + 2 h(a)} - 1}{ \beta_\omega[\omega] [\eta_p] } 
\right)^{1/p}
\end{equation*}
where $0 < \beta_\tau < 1$ and $\beta_\omega \ge 1$ are additional safety factors,
for instance $\beta_\tau = 0.75$ and $\beta_\omega = 10$. Instead of choosing $\beta_\omega$ fixed it seems worthwhile to develop a more adaptive criterion for choosing $\beta_\omega$ in the future.

Since $[\eta_p]$ and $[\omega]$ are only lower bounds and $[\tau]$ is only an upper bound for the theoretical quantities it is possible that the step size $t = [t_{\max}]$ is larger than $t_{\max}$.
Then it can happen that the Newton corrector algorithm \ref{alg:newton-corrector} rejects the initial guess since $\Theta_k > a^{2^k}$ for some $k$.
In this case a suitable step size correction formula is 
\begin{equation}\label{eq:correction-theta-0}
t' = \left(\frac{\sqrt{1 + 2h(\frac{1}{2}a)} - 1}{\sqrt{1 + 2h\left(\Theta_k^{2^{-k}} \right)} - 1} \right)^{1/p} t
\end{equation}
which is a clear reduction since $\Theta_k^{2^{-k}} > a$.

Before the path tracking algorithm is stated, the similarities and differences to the step size control developed in \cite{Telen:VanBarel:Verschelde:2019} are stated.
In \cite{Telen:VanBarel:Verschelde:2019} the authors develop an adaptive step size control similar to \eqref{eq:comp-max-step-size} with the difference that their algorithm uses instead of $[t^*]$ the step size candidate $\Delta t_1$ computed as follows.
For $\Delta t_1$ an estimate $\delta$ of the distance to the nearest path is computed based on a second-order Taylor expansion around $x(0)$.
This involves the computation of the Hessian of $H$ and multiple singular value decompositions. Then $\Delta t_1 = \left( \beta_1 \delta / [\eta_p] \right )^{1/p}$
where $\beta_1$ is a safety factor to unknown region of convergence of Newton's method, for instance $\beta_1 = 0.005$.
The cost of computing $\Delta t_1$ is $\mathcal{O}(n^4)$, whereas computing $[t^*]$ has almost no overhead.
The results in Section \ref{sec:computations} demonstrate that the approach developed in this article is much more efficient than using $\Delta t_1$ without substantially increasing the risk of path jumping.

Finally, the path tracking algorithm, Algorithm \ref{alg:path-tracking}, is stated.
It is assumed that the computational estimates $[\omega]$ and $[\mu]$ for the start solution $s$ are available.
These are either available as a result of a previous path tracking or by using Algorithm \ref{alg:model-initialization}.
\stepcounter{environment}
\begin{algorithm}[H]
\small
\caption{Path Tracking Algorithm}\label{alg:path-tracking}
\begin{algorithmic}[1]
\Require Homotopy $H(x,t): \C^n \times \C \rightarrow \C^n$, $s \in \C^n$ such that $s$ is an approximate zero of $G(x) := H(x,0)$, estimates $[\omega]$ and $[\mu]$ of $\omega$ and $\mu$, $a \in (0,\frac12]$, safety factors $\beta_\omega$ and $\beta_\tau$, minimal step size $t_{\min}$, type of Pad\'e approximant $L$.
\Ensure Approximate zero of $F(x) := H(x,1)$ or \textbf{false} if the tracking failed.
\Procedure{Track}{$H, s, [\omega], [\mu], a, \beta_\omega, \beta_\tau, t_{\min}, L$}
\State $(t, \Delta t) \gets (0, \infty)$
\State $x \gets s$
\State $\bar{u} \gets u$
\State Initialize scaling factors $D$ using $x$ and \eqref{eq:weight-init}
\While{$t < 1$}
\State Compute $x^{(1)}(t),\ldots,x^{(L+2)}(t)$ from $(x, t)$ by using the identity \eqref{eq:higher-derivatives-identity}
\State Compute $[\eta_{L+2}]$ and $[\tau]$ for $(L,1)$ Pad\'e approximant using \eqref{eq:pade-L} and \eqref{eq:errcoeff}
\State $\Delta t \gets \min 
\left(
\left[
\frac{\sqrt{1 + 2 h(a)} - 1}{ \beta_\omega [\omega] [\eta_{L+2}] } 
\right]^{1/(L+2)},
1-t,
\beta_\tau [\tau]
\right)$
\State Use $(L,1)$ Pad\'e approximant to obtain initial guess $\hat{x}$ at $t + \Delta t$ \label{line:predict}
\State Update scaling factors $D$ using $\hat{x}$ and \eqref{eq:weight-init}
\State $(\texttt{success}, \bar{x}, [\mu], [\omega], \Theta_k) \gets \Call{Newton}{\hat{x}, a, [\mu], [\omega], \bar{u}, D}$
\If{\texttt{success}}
\State $t \gets t + \Delta t$
\State $x \gets \bar{x}$
\Else
\State $\Delta t\gets \left(\frac{\sqrt{1 + 2 h(0.5 a)} - 1}{\sqrt{1 + 2 h \left(\Theta_k^{2^{-k}} \right)} - 1} \right)^{1/(L+2)} \Delta t$
\If{$\Delta t < t_{\min}$}
\State \textbf{return} \textbf{false} 
\EndIf
\State \textbf{go to} Line \ref{line:predict}
\EndIf
\If{$\bar{u} = u$ \textbf{and} $[\omega] [\mu] > a^5 h(a)$} \Comment{Use extended precision, see \eqref{eq:newton-min-acc}}
\State $\bar{u} \gets u^2$
\ElsIf{$\bar{u} = u^2$}
\State Compute estimate $\mu_u$ of $\mu(x, u, u)$ as described in Subsection \ref{subsec:comp-newton}
\If{$[\omega] \mu_u < a^7 h(a)$}
\State $\bar{u} \gets u$
\EndIf
\EndIf
\EndWhile
\State \textbf{return} $x$
\EndProcedure
\end{algorithmic}
\end{algorithm}

\section{Computational Experiments}\label{sec:computations}

In this section, numerical experiments are shown to illustrate the effectiveness of the proposed path tracking algorithm.
A prototype implementation of the path tracking algorithm together with all the data necessary to run the experiments is available at
\begin{center}
    \url{https://doi.org/10.5281/zenodo.3667414} \,.
\end{center}
The path tracking algorithm will also be implemented in version 2.0 of the Julia package \texttt{HomotopyContinuation.jl} \cite{HomotopyContinuation.jl}.

In the experiments, the implementation is compared with the state of the art. For the different solvers the following notations are used in the experiments:
\begin{center}
\begin{tabular}{ll}
\texttt{HC.jl} & The proposed algorithm, soon implemented in \texttt{HomotopyContinuation.jl}, \\
\texttt{HC.jl 1.4} & Version 1.4.0 of \texttt{HomotopyContinuation.jl}, \\
\texttt{Bertini DP} & Bertini v1.6 using double precision arithmetic (MPTYPE = 0) \cite{Bertini}, \\
\texttt{Bertini AP} & Bertini v1.6 using adaptive precision (MPTYPE = 2) \cite{Bates:Hauenstein:Sommese:Wampler:2008}, \\
\texttt{phc -u} & PHCpack \cite{Verschelde:PHCpack} v2.4.72 via \texttt{phc -u} \cite{Telen:VanBarel:Verschelde:2019} \,. \\
\end{tabular}
\end{center}
\noindent All solvers are used intentionally with the default settings unless otherwise mentioned since for a non-expert user it is very hard to understand which parameters need to be changed.
Note that Bertini uses by default a path tracking tolerance of $10^{-5}$ before and $10^{-6}$ during the endgame.
The experiments are performed on a 24 GB RAM machine with an Intel Core i5-7500 CPU working at 3.40 GHz.
All solvers use only one core for all the experiments unless stated otherwise.
The experiments are designed such that the tracked solution paths are all smooth. Therefore endgame algorithms are not necessary.
In all experiments the implementation of the proposed algorithm uses a $(2,1)$ Pad\'e approximant with parameters $a=0.2$, $\beta_\omega = 10$ and $\beta_\tau = 0.75$.

\begin{figure}[h]
\label{fig:hyperbolas}
\includegraphics[width=0.5\textwidth]{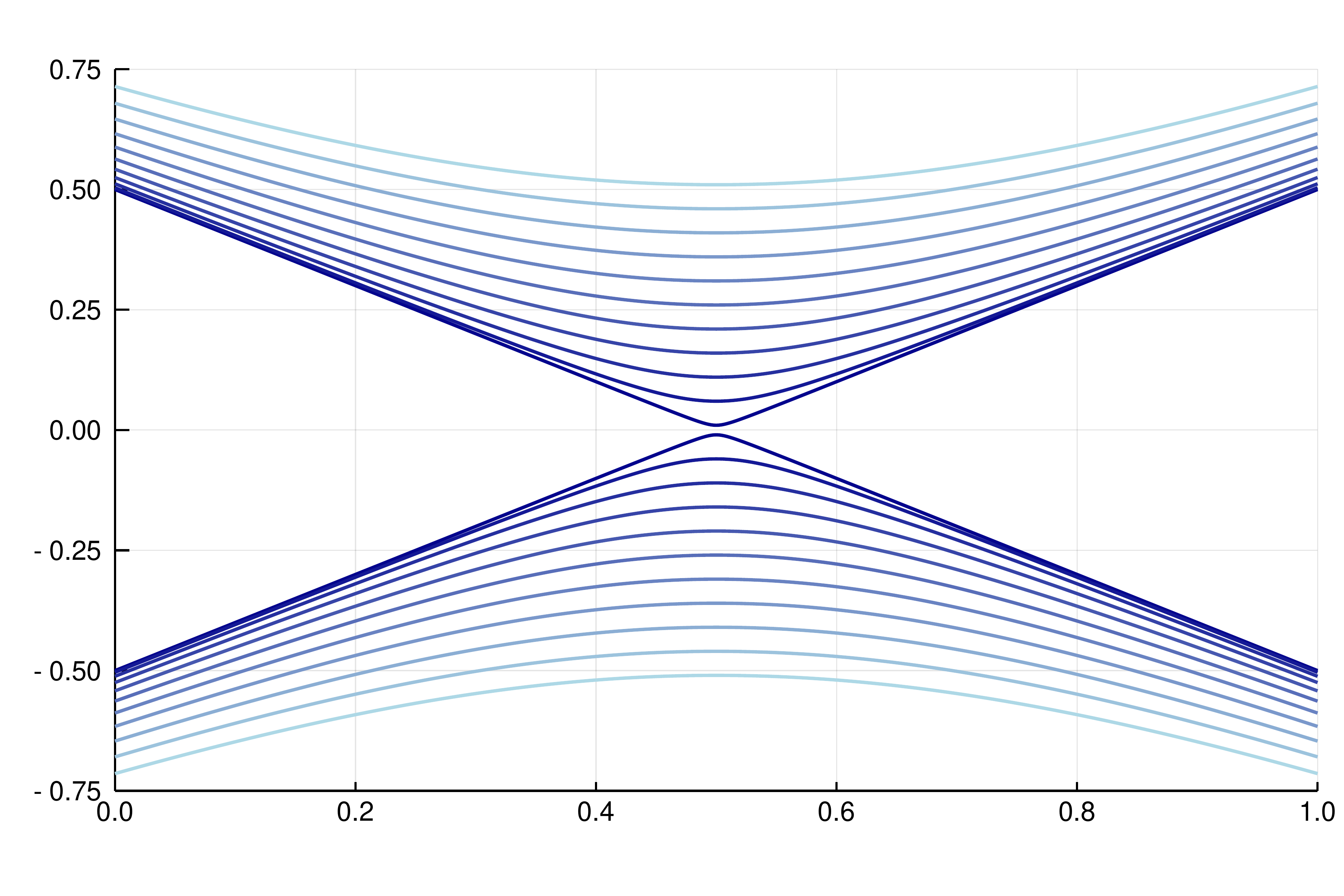}
\caption{Family of hyperbolas defined by \eqref{eq:hyperbola-family} for $\rho$ between $0.51$ and $0.01$.}
\end{figure}
\subsection{A Family of Hyperbolas}\label{subsec:hyperbolas} 
This is an experiment from \cite{Telen:VanBarel:Verschelde:2019} which constructs a situation which causes danger for path jumping.
Consider a family of hyperbolas parametrized by a real parameter $\rho$ and defined by
\begin{equation}\label{eq:hyperbola-family}
    H_\rho(x,t) = x^2 - (t- 1 / 2)^2 - \rho^2 = 0 \,.
\end{equation}
See Figure \ref{fig:hyperbolas} for an illustration.
\noindent For general $t^*$ the homotopy $H_\rho$ has the two solutions
\begin{equation*} 
    \pm \sqrt{\large(t^* - 1/ 2 \large)^2 + \rho^2} \,.
\end{equation*}
Note that the Jacobian $H_x = \frac{\partial}{\partial x} H$ is singular for $t^* = \frac12 + \rho\sqrt{-1}$ and $x=\sqrt{(t^* - \frac12)^2 + \rho^2}$.
Consider $\rho \neq 0$ here, such that $x(t): [0,1] \rightarrow \C$ is a smooth path.
The smaller $|\rho|$, the closer the two branch points move to the line segment $[0,1]$. Figure \ref{fig:hyperbolas} shows that as the value of $\rho>0$ decreases, the two solution paths approach each other for parameter values $t^* \approx 0.5$ which causes danger for path jumping.
The experiments confirm this with the results shown in Table \ref{tab:hyperbolas}.
Note that \texttt{HC.jl 1.4} and \texttt{Bertini AP/DP} are prone to path jumping. The new path tracking algorithm, on the other hand, handles this situation correctly.

\begin{table}[h]
\centering
\begin{tabular}{l|lllllll}
\diagbox{Solver}{$k$} & 1 & 2 & 3 & 4 & 5 & 6 & 7 \\ \hline
\texttt{Bertini DP} & \cmark & \cmark & \cmark & \xmark & \xmark & \xmark & \xmark \\
\texttt{Bertini AP} & \cmark & \cmark & \cmark & \cmark & \cmark & \xmark & \xmark \\
\texttt{phc -u} & \cmark & \cmark & \cmark & \cmark & \cmark & \cmark & \cmark \\
\texttt{HC.jl 1.4} & \cmark &\xmark & \xmark & \xmark & \xmark & \xmark & \xmark \\
\texttt{HC.jl} & \cmark & \cmark & \cmark & \cmark & \cmark & \cmark & \cmark \\
\end{tabular}
\caption{Results of the experiment of Subsection \ref{subsec:hyperbolas} for $\rho = 10^{-k}$ and $k = 1, \ldots, 7$. A `\xmark' indicates that path jumping happened.}
\label{tab:hyperbolas}
\end{table}

\subsection{Alt's problem}\label{subsec:fourbar}

Alt's problem, formulated in 1923, is to count the number of four-bar linkages whose coupler curve interpolates nine general points in the plane.
In 1992, Morgan, Sommese, and Wampler \cite{Wampler:Morgan:Sommese:1992} provided a numerical proof to Alt's problem that there are generically 1442 non-degenerate four-bar linkages. Due to Roberts cognates and a two-fold symmetry, the resulting polynomial system generically has 8652 regular solutions. Here the formulation of the polynomial system as an affine polynomial system in 24 variables and the 16 parameters $(\delta, \hat{\delta}) \in \C^8 \times \C^8$ is considered.
Since the problem is formulated in isotropic coordinates, the physically meaningful configurations correspond to choices of parameters such that $\delta$ and $\hat{\delta}$ are complex conjugates.
Consider the `general' situation where solutions are tracked from generic parameter values $\delta_1 \in \C^8 \times \C^8$ to generic physically meaningful parameter values $(\delta_0, \bar{\delta_0})$ with $\delta_0 \in \C^8$.
The results in Table \ref{tab:fourbar_generic} show that even this general situation results in numerically challenging paths which \texttt{Bertini AP} cannot handle with its default settings.
After decreasing the path tracking tolerance to $10^{-8}$, in half of the cases all solutions are found. The proposed algorithm, on the other hand, reliably computes all 1442 solutions without any path jumping in a fraction of the time.

\begin{table}[h]\centering 
\footnotesize
\begin{tabular}{lcrrrrrrrr}\toprule
& & \multicolumn{4}{c}{runtime (seconds)} & \phantom{abc} & \multicolumn{3}{c}{\# solutions} \\
\cmidrule{3-6} \cmidrule{8-10}
& \texttt{tol} & mean & median & min & max &  & median & min & max \\ \midrule
\texttt{HC.jl} & - & 6.38 & 6.41 & 5.08 & 7.12 & &  1442 & 1442 & 1442 \\
\texttt{Bertini AP} & default & 1931.63 & 1832.59 & 1120.43 & 2695.13 & & 1436 & 1433 & 1441 \\
\texttt{Bertini AP} & 1e-8 & 10950.09 & 11218.97 & 8371.76 & 12115.66 & &  1441 & 1437 & 1442 \\
\bottomrule
\end{tabular}
\caption{Results for 10 runs of Alt's problem using the same generic start solutions to a generic physically meaningful configuration. The \texttt{tol} column refers to the assigned path tracking tolerance.}
\label{tab:fourbar_generic}
\end{table}

\begin{figure}[b]
\includegraphics[width=0.7\textwidth]{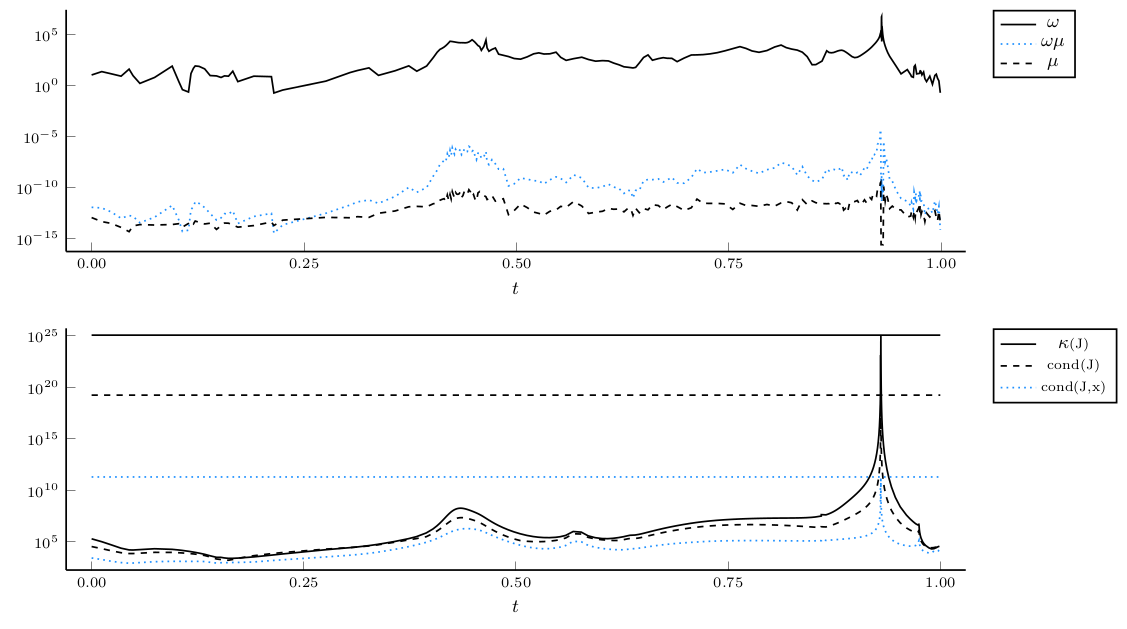}
\caption{Behavior of the path tracking algorithm along a single challenging path.}
\label{fig:fourbar-bad-path}
\end{figure}
Figure \ref{fig:fourbar-bad-path} illustrates the behavior of the proposed algorithm for one particular numerically challenging path.
This path closely passes at around $t = 0.93$ a singularity.
As expected, this increases the Lipschitz constant $\omega$ and the limit accuracy $\mu$.
The limit accuracy decreases sharply as soon as the algorithm switches to extended precision.
After passing the problematic region, the algorithm quickly switches back to only using double precision arithmetic, as indicated by a sharp increase of $\mu$.
The lower part of Figure $\ref{fig:fourbar-bad-path}$ depicts different values associated with the Jacobian $J = H_x(x(t),t)$ along the path:
the condition number $\kappa(J)$, the componentwise relative condition number $\mathrm{cond}(J)$ and $\mathrm{cond}(J,x)$.
See \cite[Sec. 7]{Higham:2002} for the definitions.
The values of $\kappa(J)$ and $\mathrm{cond}(J)$ obtain a maximum of $2.1\times 10 ^{25}$ resp. $1.6\times 10 ^{19}$.
From these values, it seems hopeless to track the path in double precision arithmetic due to the well-known rule of thumb that one expects to lose around $\log_{10}(\kappa(J))$ digits of accuracy in the linear system solving.
However, the componentwise relative forward error of the computed Newton updates is only governed by the much tamer $\mathrm{cond}(J,x)$ which is at most $5.8 \times 10^{10}$.
This explains why it is still possible to track the path by only using mixed precision iterative refinement.
Using the proposed path tracking algorithm the path needs $253$ steps in total with only two rejected steps and a total runtime of 13 milliseconds.
Trying to compute the path with \texttt{Bertini AP} results in a path failure after around 90 seconds and over 5000 steps due to supposedly insufficient precision of at most 1024 bits.
After changing the required tracking tolerance to $10^{-12}$ \texttt{Bertini AP} successfully tracks the path in $120$ seconds.

\subsection{Steiner's Conic Problem}

A classic problem in enumerative geometry is Steiner's conic problem.
It asks: How many plane conics are tangent to five given conics in general position?
See \cite{Eisenbud:Harris:2016} for historic remarks and a modern intersection theory treatment of this classic problem.
Here the formulation of Steiner's conic problem from \cite{Breiding:Sturmfels:Timme:2020} is used where the resulting polynomial system has 15 variables and 30 parameters.
To test the path tracking algorithm consider the case of a parameter homotopy \cite{Morgan:Sommese:89} from generic complex parameters $c \in \C^{30}$ to generic real parameters $r \in \R^{30}$.
The results for 50 parameter homotopies are shown in Table \ref{tab:steiner}.
As for Alt's problem, the proposed path tracking algorithm handles all instances without any failure or path jumping.
However, even these generic instances pose problems for other path tracking algorithms with \texttt{Bertini AP} losing solutions almost always solutions using the default settings. After the path tracking tolerance is manually decreased to $10^{-8}$ in almost all instances all solutions are found.

\begin{table}[h]\centering 
\footnotesize
\begin{tabular}{lcrrrcrrrr}\toprule
& & \multicolumn{4}{c}{runtime (seconds)} & \phantom{abc} & \multicolumn{3}{c}{\# solutions} \\
\cmidrule{3-6} \cmidrule{8-10}
& \texttt{tol} & mean & median & min & max &  & median & min & max \\ \midrule
\texttt{HC.jl} & & 2.47 & 2.49 & 1.64 & 3.11 & & 3264 & 3264 & 3264 \\
\texttt{Bertini DP} & default & 34.42 & 34.57 & 23.85 & 45.34 & &  3189 & 3094 & 3216 \\
\texttt{Bertini AP} & default & 130.53 & 126.61 & 78.70 & 251.01 & & 3261 & 3256 & 3264 \\
\texttt{Bertini AP} & 1e-8 & 688.50 & 691.59 & 368.36 & 1125.72 & & 3264 & 3261 & 3264 \\
\bottomrule
\end{tabular}
\caption{Results for 50 runs of Steiner's problem from generic complex parameter values to generic real parameter values. The \texttt{tol} column refers to the assigned path tracking tolerance.}
\label{tab:steiner}
\end{table}

\subsection{The Katsura Family}
The {\tt katsura} family of systems is named after the problem
posed by Katsura~\cite{Katsura:94}, see~\cite{Katsura:90} for a description
of its relevance to applications.
The {\tt katsura}-$n$ problem consists of $n$ quadratic equations
and one linear equation. 
The number of solutions equals $2^n$, the product of the degrees of all polynomials in the system.

\setlength{\tabcolsep}{4.0pt}
\begin{table}[h]
\footnotesize
\begin{center}
\begin{tabular}{c|r|rr|r|r|r|r}
$n$ & \multicolumn{1}{c|}{\#sols} & \#real
& \multicolumn{1}{c|}{\#imag}
& \multicolumn{2}{c|}{\texttt{HC.jl} time (seconds)}
& \multicolumn{2}{c}{\texttt{phc -u} time (seconds)} \\ \hline
12 & 4,096 & 582 & 3,514 & 2.85 & 3s & 7.925E+01 & 1m 19s \\
13 & 8,192 & 900 & 7,292 & 7.22 & 8s & 2.081E+02 & 3m 28s \\
14 & 16,384 & 1,606 & 14,778 & 26.96 & 19s & 5.065E+02 & 8m 27s \\ 
15 & 32,768 & 2,542 & 30,226 & 42.87 & 45s & 1.456E+03 & 24m 16s \\
16 & 65,536 & 4,440 & 61,096 & 113.21 & 1m 59s &4.156E+03 & 1h \phantom{0}9m 16s \\
17 & 131,072 & 7,116 & 123,956 & 231.49 & 4m 48s &1.001E+04 & 2h 46m 50s \\ 
18 & 262,144 & 12,458 & 249,686 & 855.90 & 14m 16s & 2.308E+04 & 6h 24m 15s \\
19 & 524,288 & 20,210 & 504,078 & 1330.31 & 22m 10s & 5.696E+04 & 15h 49m 20s \\
20 & 1,048,576 & 35,206 & 1,013,370 & 3772.08 & 1h \phantom{0}2m 52s & 1.317E+05 & 36h 34m 11s
\end{tabular}

\caption{Runtime of \texttt{HC.jl} on a \emph{single} core and wall clock time on 44 cores of \texttt{phc -u} for the {\tt katsura} problem.}
\label{tabbenchkatsura}
\end{center}
\end{table}

Table~\ref{tabbenchkatsura} summarizes the characteristics
and run times on {\tt katsura}-$n$, for $n$ ranging from 12 to 20 on a \emph{single} CPU.
For reference the runtime of \texttt{phc -u} is stated as reported in \cite{Telen:VanBarel:Verschelde:2019}.
Note that these computations were performed on two 22-core 2.2 GHz processors, i.e., using 44 CPUs in total.
In \cite{Telen:VanBarel:Verschelde:2019} it is also mentioned that the computation required the use of homogeneous coordinates since otherwise intermediate coordinate growth resulted in tracking failures for some paths.
However, this is not necessary if a dynamic rescaling of the variables is performed as derived in Section \ref{subsec:comp-newton}.

\section{Conclusion and Future Work}
This article proposed a mixed precision path tracking algorithm for numerical path tracking in polynomial homotopy continuation.
The results of the computational experiments demonstrate that the algorithm is very robust, can handle numerically challenging situations and, at the same time, is faster than existing software implementations.
An implementation is available, and the algorithm will also be integrated into version 2.0 of \texttt{HomotopyContinuation.jl}.
It is expected that the techniques in this article are also helpful in developing more efficient and robust endgame algorithms to deal with singular solutions and diverging paths.

\section*{Acknowledgement}
The author wants to thank Paul Breiding, Peter Deuflhard, Michael Joswig, Simon Telen, and Marc Van Barel for helpful discussions.
The author was supported by the Deutsche Forschungsgemeinschaft (German Research Foundation) Graduiertenkolleg {\em Facets of Complexity} (GRK~2434).

\bibliographystyle{alpha} 
\bibliography{bibliography.bib}

\end{document}